\documentclass[12pt,oneside,a4paper,leqno]{article}
\usepackage[utf8]{inputenc}
\usepackage[english]{babel}
\usepackage[all]{xy}

\CompileMatrices
\usepackage{verbatim}

\newdir{ >}{!/10pt/@{ }*@{>}}
\newdir{< }{!/10pt/@{ }*@{<}}

\usepackage{graphicx}
\usepackage{amsfonts,amssymb}
\usepackage[centertags]{amsmath}

\usepackage{amsthm} %
\newcounter{theorems}

\theoremstyle{plain}

\swapnumbers
\newcounter{lemma}
\numberwithin{equation}{section}

\newtheoremstyle{par}%
     {\topsep}%
     {\topsep}%
     {\itshape}%
     {}%
     {\bfseries}%
     {}%
     {.5em}%
     {}%

\newtheoremstyle{parrm}%
     {\topsep}%
     {\topsep}%
     {\normalfont}%
     {}%
     {\itshape}%
     {}%
     {.5em}%
     {}%

\theoremstyle{plain}
\numberwithin{equation}{section}
\newtheorem{lemma}[equation]{Lemma}
\newtheorem{theo}[equation]{Theorem}

\theoremstyle{definition}

\newtheorem{example}[equation]{Example}

\theoremstyle{remark}
\newtheorem{remark}[equation]{Remark}

\theoremstyle{par}

\newtheorem{propo}[equation]{}

\theoremstyle{parrm}

\makeatletter
\def\tagform@#1{\maketag@@@{\ignorespaces#1\unskip\@@italiccorr}}

\makeatother

\usepackage{paralist}
\usepackage{comment}
\usepackage{subcaption} %
\newcommand{\RR}{\mathbb{R}}

\newcommand{\PP}{\mathbb{P}}

\newcommand{\from}{\colon}

\usepackage{bm}
\newcommand{\simbolovettore}[1]{{\boldsymbol{#1}}}

\newcommand{\vm}{\simbolovettore{m}}

\newcommand{\vq}{\simbolovettore{q}}

\newcommand{\vv}{\simbolovettore{v}}
\newcommand{\vw}{\simbolovettore{w}}
\newcommand{\vx}{\simbolovettore{x}}

\newcommand{\vQ}{\simbolovettore{Q}}

\newcommand{\vL}{\simbolovettore{L}}
\newcommand{\vY}{\simbolovettore{Y}}

\newcommand{\zero}{\boldsymbol{0}}

\newcommand{\Fix}{\operatorname{Fix}}

\usepackage{tikz}
\usetikzlibrary{matrix}

\newcommand{\norm}[1]{\lVert{#1}\rVert}

\newcommand{\conf}[2]{\mathbb{F}_{#1}(#2)}

\newcommand{\CH}{\operatorname{CH}}

\begin{document}
\pagenumbering{arabic}

\title{%
Fixed points and the inverse problem for central configurations
}

\author{D.L.~Ferrario
}

\date{%
\today}
\maketitle

\begin{abstract}
Central configurations play an important role in the dynamics of the $n$-body
problem: they occur as relative equilibria and as asymptotic configurations in
colliding trajectories.  We illustrate how they can be found as projective
fixed points of self-maps defined on the shape space, and some results on the
inverse problem in dimension $1$, i.e. finding (positive or real) masses which
make a given collinear configuration central.
This survey article introduces readers to the recent results of the author,
also unpublished, showing an application of the fixed point theory.

\end{abstract}

\section{Introduction}

Let $n\geq 3$ be and integer, and $m_1,\ldots, m_n$ positive parameters,
masses. Given a dimension $d\geq 2$, a configuration of $n$ points 
in $\RR^d$ is a $n$-tuple $\vq=(\vq_1,\ldots, \vq_n)$ with $\vq_j\in \RR^d$
for all $j$ and $\vq_i \neq \vq_j$ whenever $i\neq j$. 
Spaces of configurations (configuration spaces) have been the object
of much of study in recent decades (see Fadell--Husseini 
\cite{fadellGeometryTopologyConfiguration2001}
for a topological point of view and some deep and interesting consequences). The set of all configurations is 
  denoted $\conf{n}{\RR^d}$, following Fadell--Husseini notation. 
Topological properties of configuration spaces have consequences on the study 
of dynamical systems of $n$ point particles interacting with
mutual forces (the $n$-body problem: see for example 
\cite{SmaleTopologymechanicsII1970}).
One of the most immediate problem occurring in this context is the problem
of finding and classifying \emph{central configurations}: 
if the particles interact under a potential $U$
  defined, for a given homogeneity parameter $\alpha>0$, as 
\[
U = \sum_{i<j} \dfrac{m_im_j}{\norm{\vq_i - \vq_j}^{\alpha}},
\]
then central configurations are configurations of points such that there exists 
a (negative) constant $\lambda$ such that 
\begin{equation}\label{eq:CC}
\lambda m_i \vq_i = - \alpha \sum_{j\neq i} \dfrac{m_i m_j}{\norm{\vq_i - \vq_j}^{\alpha+2}} (\vq_i -\vq_j) 
= \dfrac{\partial U}{\partial \vq_i}
\end{equation}
for $i=1,\ldots, n$. 
For $\alpha=1$, it is the Newtonian gravitational interaction.
There have been an active line of research on central configurations since decades: see for example 
\cite{Moultonstraightlinesolutions1910},
\cite{Buchanancertaindeterminantsconnected1909}
\cite{SmaleTopologymechanicsII1970}),
\cite{xiaCentralConfigurationsMany1991},
\cite{moeckelCentralConfigurations1990},
\cite{hamptonFinitenessRelativeEquilibria2006},
\cite{AlbouyFinitenesscentralconfigurations2012},
\cite{moeckelCentralConfigurations2015},
\cite{AlbouyOpenProblemAre2015}. 
The purpose of this article is to survey some recent results (mostly of the
author) in this field, with applications of fixed point theory.

\section{Central configurations as fixed points: self-maps and multi-valued self maps}

Since the potential $U$ and its gradient $\nabla U$ are homogeneous in $\vq$, 
equation \eqref{eq:CC} can be re-written as a fixed point problem 
defined on the unit ellipsoid
$S = \{ \vq \in \conf{n}{\RR^d} : \sum_{i=1}^n m_i \norm{\vq_i}^2 = 1 \}$: 
\begin{equation}\label{eq:CC2}
-\frac{\lambda}{\alpha} \vq_i = \sum_{j\neq i} \dfrac{m_j}{\norm{\vq_i - \vq_j}^{\alpha+2}} (\vq_i -\vq_j) 
= -\frac{1}{\alpha m_i} \dfrac{\partial U}{\partial \vq_i}
\end{equation}
If $\nabla_M$ denotes the gradient with respect to the mass-metric $\langle -,-\rangle_M$ 
defined on the tangent vectors of $\conf{n}{\RR^d}$ as 
\[
\langle \vv, \vw \rangle_M = \sum_{i=1}^n m_i \vv_i \cdot \vw_i,
\]
where $\vv_i \cdot \vw_i$ is the standard euclidean scalar product in $\RR^d$, 
then equation \eqref{eq:CC2} can be written as 
\begin{equation}
\label{eq:CC3}
\vq = \frac{-\nabla_M U}{ \norm{\nabla_M U}_M } = F(\vq),
\end{equation}
for a map $F \from \conf{n}{\RR^d} \to \RR^{nd}$. 
  Here $\norm{-}_M$ is the norm associated to the mass-metric scalar product $\langle -,-\rangle_M$,
hence  $S$ is the unit sphere in $\conf{n}{\RR^d}$ with respect to the norm mass-metric norm. 
This notation is standard, and highlights the fact that here masses are a chosen (and fixed) parameter 
of the problem. For the inverse problem, below, this notation will be not necessary since
masses will be part of the solution. 

Now consider the euclidean group of symmetries of $\RR^d$: it acts diagonally (component-by-component) 
on the
configuration space $\conf{n}{\RR^d}$. Since $U$ is invariant with respect to translations in $\RR^d$,
the map $F$ is invariant with respect to translations: $F(\vq) = F(\vq + \vv)$, 
whenever $\vv$ is a vector of type $\vv_1=\vv_2=\ldots = \vv_n$. Hence, the image
$F(\vq)$ is orthogonal to any such a diagonal $\vv$ (with respect to $\langle -,-\rangle_M$),
i.e. it belongs to the subspace
\[
X_0 = \{ \vv \in (\RR^d)^n : \sum_{i=1}^n m_i \vv_i  = \zero \}. 
\]
Since any fixed point (central configuration) must belong to $X_0$, 
we can restrict the inertia ellipsoid to $X_0$, and consider the restricted map
$F_0 \from S_0 = S \cap X_0 \to \overline S_0$. 
Here we must take the closure $\overline S_0$ because it is not guaranteed that $F(\vq)$ 
is collision-free for every $\vq\in \conf{n}{\RR^d}$. 

Now, the potential $U$ is $SO(d)$-invariant, 
where again the $SO(d)$ action is diagonal
on $\conf{n}{\RR^d}$, 
and hence its $M$-gradient and the function $F$ are $SO(d)$-equivariant:
for each $g\in SO(d)$ one has $F(g\vq) = g F(\vq)$. 
This means that, if $\pi \from X_0 \to X_0/SO(d)$ denotes the projection (and the same 
for $S_0 \to S_0/SO(d)$), 
the map $F$ induces a map $f$ on the quotient $S_0/SO(d)$ 
\begin{equation}
\xymatrix{
  S_0 \ar[r]^F \ar[d]^\pi & \overline S_0 \subset X_0 \ar[d]^\pi  \\
S_0/SO(d) \ar[r]^f & \overline S_0/SO(d). %
}
\end{equation}

The results in following proposition were proved in 
\cite{ferrarioPlanarCentralConfigurations2007},
\cite{FerrarioFixedpointindices2015},
\cite{ferrarioCentralConfigurationsMorse2017}.
See also \cite{ferrarioCentralConfigurationsMutual2017}), 
where a combinatorially cohomological approach on coordinates
was used to simplify the mass-metric projections, which we are going to use 
later below in \eqref{eq:Qij}.
\begin{propo}
The map $f$ is well-defined, compactly fixed, and 
\[
\pi \Fix(F)  = \Fix f.
\] 
If $[\vq]=\pi(\vq)$ is an isolated fixed-point of $f$, with  maximal isotropy stratum 
of $S_0/SO(d)$, then  its fixed point index is $(-1)^\mu$,
where $\mu$ is the Morse index of $U$ at $\vq$. 
\end{propo}

\begin{remark}
\label{rem:Qij}
Consider the following variables, for all $i,j=1,\ldots, n$
\begin{equation}
\label{eq:Qij}
\vq_{ij} = \vq_i - \vq_j ; \quad
\vQ_{ij} = \begin{cases}
\frac{\vq_{ij}}{\norm{\vq_{ij}}^{\alpha+2} } & \text{ if } i\neq j \\
0 & \text{ if } i=j.
\end{cases} 
\end{equation}
Then 
\[
F(\vq)
= \dfrac{\sum_{j=1}^n m_j \vQ_{ij} }{\norm{ \sum_{j=1}^n m_j \vQ_{ij} }_M }.
\]
Hence $F(\vq)$ belongs to the positive cone generated by the column
vectors if the skew-symmetric matrix with entries $\vQ_{ij}$. 
\end{remark}

\begin{figure}
\centering
    \begin{subfigure}[t]{0.46\textwidth}
        \centering
        \includegraphics[width=0.8\textwidth]{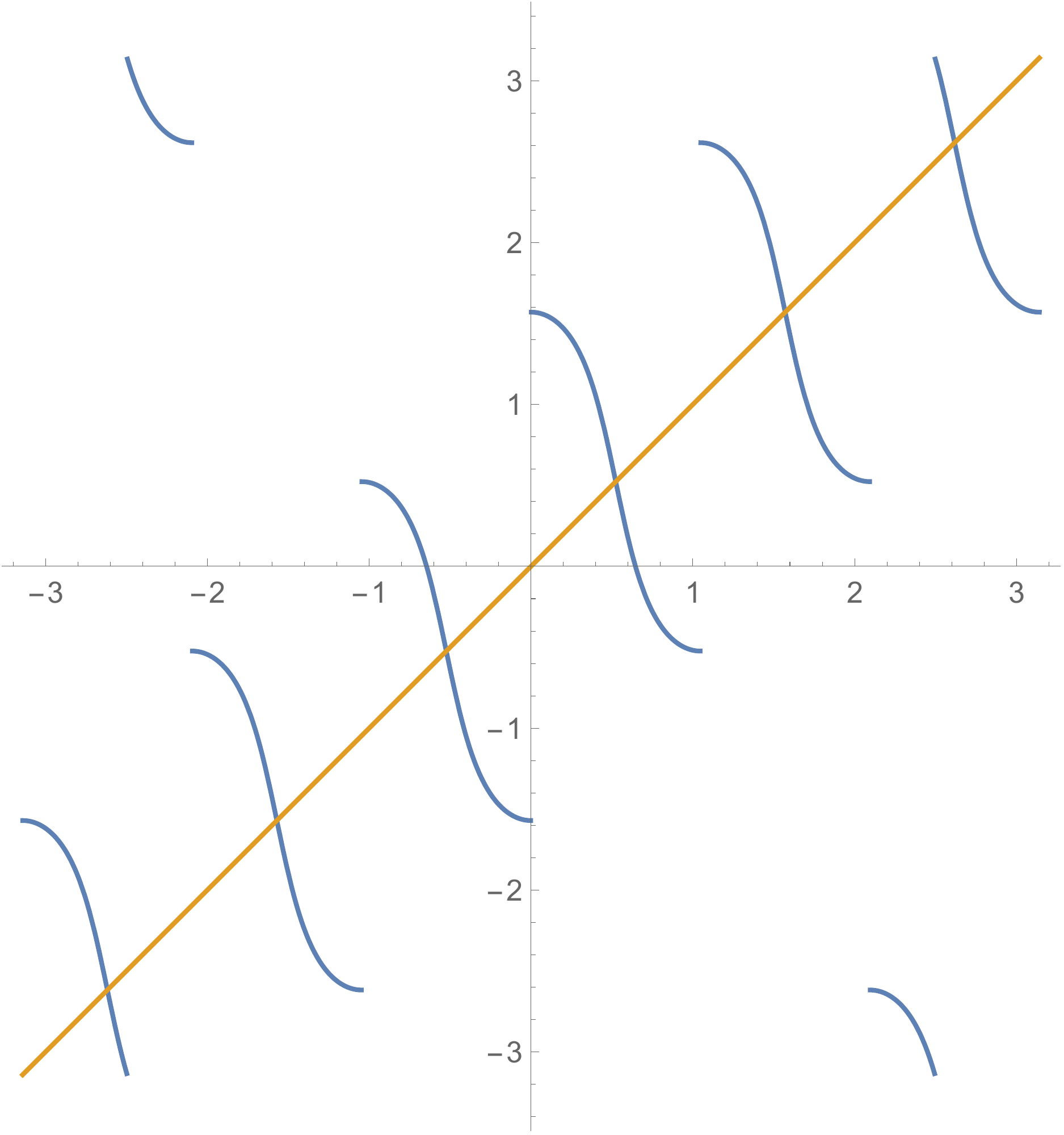}
        \caption{The discontinuous function $f\from S_0 \to S_0$ for $n=3$ and $m_i=1$}
    \end{subfigure}%
    ~
    \begin{subfigure}[t]{0.46\textwidth}
        \centering
        \includegraphics[width=0.8\textwidth]{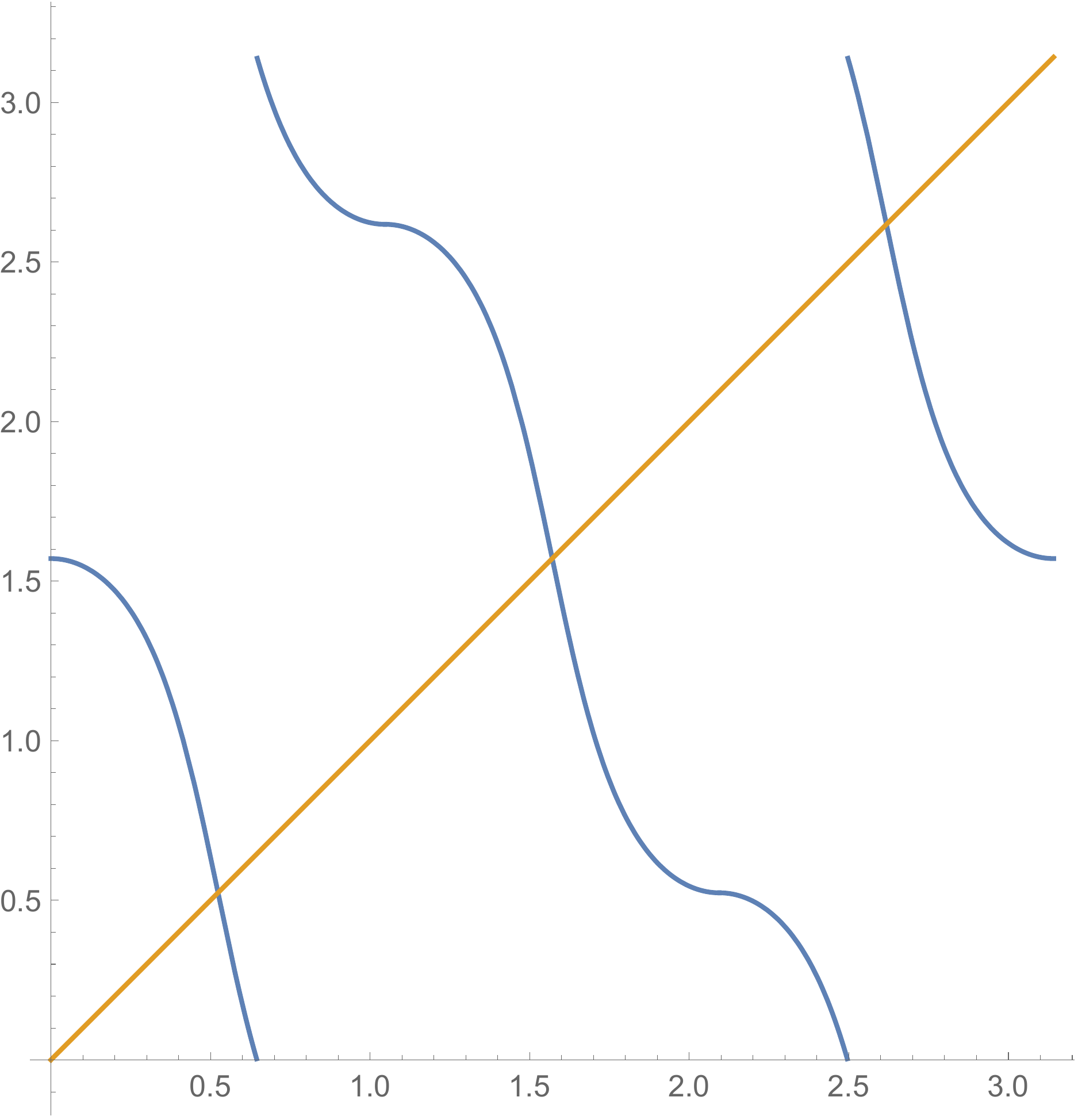}
        \caption{The continuous projective function $\bar f \from \PP^1(\RR) \to \PP^1(\RR)$ for $n=3$ and $m_i=1$}
    \end{subfigure}%
\caption{The two functions $f$ and $\bar f$ for $n=3$: on the left $\mod 2\pi$, on the right $\mod \pi$}
\label{fig:1}
\end{figure}

\begin{example}
If $d=1$, then $SO(d)$ is trivial, and the dimension of $S_0$ is $n-2$: 
it has $n!$ connected components, one for each strict ordering of the coordinates 
$q_1,\ldots, q_n$. If $n=3$, then the map $f$ can not  be extended by continuity on 
$\overline S_0$, so that to be a genuine self-map $S^1 \to S^1$. In fact,
the matrix $\vQ_{ij}$ is a $3\times 3$ skew-symmetric matrix 
with real entries and $F(\vq)$ is the renormalization of the vector
obtained by matrix-vector multiplication
\[
\begin{bmatrix}
0 & Q_{12} & Q_{13} \\
Q_{21} & 0 & Q_{23} \\
Q_{31} & Q_{32} & 0 
\end{bmatrix}
\begin{bmatrix}
m_1 \\ m_2 \\ m_3 
\end{bmatrix}.
\]
Consider the component $q_1>q_2>q_3$, which means $q_{12}>0$, $q_{23}>0$. 
If $q_{12} \to 0^+$, then $F(\vq)$ tends to the renormalization of the vector
(divide by $Q_{12}$ which is positive and goes to $+\infty$ as $q_{12}\to 0^+$)
\[
\begin{bmatrix}
0 & 1  & 0  \\
-1  & 0 & 0  \\
0  & 0 & 0 
\end{bmatrix}
\begin{bmatrix}
m_1 \\ m_2 \\ m_3 
\end{bmatrix}.
\]
On the other, if we approach the collision $q_1=q_2$ from the 
component $q_2>q_1>q_3$, we must divide by $Q_{21}$ which is positive and goes to $+\infty$ 
as $q_{21} \to 0^+$, the limit is 
\[
\begin{bmatrix}
0 & -1  & 0  \\
1  & 0 & 0  \\
0 & 0 & 0 
\end{bmatrix}
\begin{bmatrix}
m_1 \\ m_2 \\ m_3 
\end{bmatrix}.
\]
In order to define a map we need to consider the antipodal map $a \from S_0 \to S_0$,
the corresponding group action, and the quotient  $S_0/\pm$. 
It is easy to see that $f$ is $\pm$-equivariant, and it induces 
a map 
\[
\bar f \from S_0/\pm \subset \PP^{n-2}(\RR) \to \overline S_0/\pm.
\]
This map now can be extended to a continuous map $S_0/\pm \to S_0/\pm$. 
For $n=1$, it is the map $\bar f\from S^1 \to S^1$ represented in figure 
\ref{fig:1}.
The same happens for any $n>3$: $S_0$ is an open subspace of the sphere of dimension $n-2$ (with $n!$
components), which projects
onto an open subspace of the projective space $\PP^{n-2}(\RR)$ (with $\frac{n!}{2}$ components). 
\end{example}

\begin{remark}
Fixed point indices sum up to form the (local) Lefschetz number, and critical point Morse
indices yield the Morse polynomial. Both give topological estimates on the number 
of central configurations, provided they are non-degenerate (and here
 Morse non-degenerate on the quotient, which is equivalent to say that 
 the Jacobian of $f$ is non-degenerate at fixed points). 
The computations involve homology computations on configuration spaces
or on projective configuration spaces.  
See \cite{PalmoreClassifyingrelativeequilibria1973},
\cite{PalmoreClassifyingrelativeequilibria1975}, 
\cite{PacellaCentralconfigurationsNbody1987},
\cite{McCordPlanarcentralconfiguration1996},
and \cite{MerkelMorseTheoryCentral2008}.
\end{remark}

Now, for each choice of masses there is a (compact) subset of central configurations in $S_0$. 
It can be proven that the set is non-empty (the minima of the potential yield 
central configurations). On the other hand, as already Moulton set forth,
there is the \emph{inverse} problem of central configurations: 
We consider, after Albouy--Moeckel \cite{AlbouyInverseProblemCollinear2000}, 
the inverse problem for central
configurations: given a configuration of $n$ bodies, find positive
masses which make it central. 
We consider only collinear configurations ($d=1$), and we refer for details and 
further references on the inverse problem to 
\cite{AlbouyInverseProblemCollinear2000},
\cite{OuyangCollinearCentralConfiguration2005},
\cite{Xieanalyticalproofcertain2014},
\cite{DavisInverseproblemcentral2018a},
\cite{ferrarioPfaffiansInverseProblem2020}.

\section{The inverse collinear central configurations problem: simplex-valued maps}  

Now we reformulate the inverse collinear problem in a quotient space. 
First, re-define $X_0 = \{ \vq\in \conf{n}{\RR} : \sum_{i=1}^n q_i = 0 \}$: 
now it does not depend on the masses. 
Consider the orthogonal projection of $\Pi\from\conf{n}{\RR} \to X_0$ (with respect to 
the 
standard euclidean scalar product in $\RR^n$, i.e. the projection 
parallel to the vector $\vL$). 
Then  the inverse problem for a configuration $\vq\in \conf{n}{\RR}$ has solutions (i.e. positive masses
such that the configuration is central with respect with this choice of masses and 
its center of mass) if and only if there exists $\vm \in \RR^n$ with positive coefficients
such that  $\Pi Q \vm = \Pi \vq$.

The space $X_0$ has dimension $n-1$: let $X_1$ denote the $(n-2)$-dimensional subspace  
\[
X_1 = \{ \vq \in X_0 : q_1 - q_n = 1 \}.
\]
Consider the open cone 
\[
X_0^+ = \{ \vq \in X_0 : q_1> q_2>\ldots > q_n \}.
\]
It is homeomorphic to $(0,+\infty) \times \mathring{\Delta}^{n-2}$,
where $\Delta^{n-2}$ is the closure of  $X_1\cap X_0^+ \}$ in $X_1$:
it is a $(n-2)$-simplex, and its
$n-1$ faces
are given by  the  equations
\[
q_1=q_2, \quad \ldots \quad 
q_{n-1}=q_n;
\]
the simplex $\Delta^{n-2}$ itself is given by the inequalities $q_1\geq q_2\geq \ldots \geq q_n$
with $\sum_{i=1}^n q_i = 0 $ and $q_1 -q_n = 1$. 

A simple computation gives the following lemma (cf. 
\cite{ferrarioPfaffiansInverseProblem2020}
for details and \cite{ferrarioCentralConfigurationsMutual2017} for more on mutual differences).
\begin{lemma}
Let $x_i=q_i-q_{i+1}$, for $i=1,\ldots, n-1$. Then $\vx=(x_i)$ are 
a linear system of coordinates on $X_0$, and $\vq\in X_0^+$ if and only if 
$x_i>0$ for $i=1,\ldots, n-1$.
Moreover, the simplex $\Delta^{n-2}$ is 
the standard $n-2$ simplex in $\RR^{n-1}$,
and $x_i$ are barycentric coordinates 
with respect to its vertices.  
\end{lemma}

In $\vx$-coordinates then for a configuration $\vq\in X_0$ 
there is a solution to the inverse problem if and only if there exists
$\vm \in \RR^n$, all positive, such that 
\begin{equation}
\label{eq:matrixY}
\begin{bmatrix}
x_1\\
x_2\\
\vdots \\
x_{n-1}
\end{bmatrix}
= 
\begin{bmatrix}
Q_{11} - Q_{21} & Q_{12} - Q_{22} & \ldots & Q_{1n}-Q_{2n} \\
Q_{21} - Q_{31} & Q_{22} - Q_{32} & \ldots & Q_{2n}-Q_{3n} \\
\vdots & \vdots &  \ddots & \vdots \\
Q_{n-1,1} - Q_{n,1} & Q_{n-1,2} - Q_{n,2} & \ldots & Q_{n-1,n}-Q_{n,n} \\
\end{bmatrix}
\begin{bmatrix}
m_1\\m_2\\\vdots\\m_n
\end{bmatrix}
\end{equation}
Let $Y$ denote the matrix with entries $Q_{i,j}-Q_{i+1,j}$.
Note that the sums of the entries in the $j$-th column is 
\begin{equation}
\label{eq:unknown}
\sum_{i=1}^{n-1} Q_{ij} - Q_{i+1,j} = Q_{1,j} - Q_{n,j} = Q_{1,j} + Q_{j,n} > 0.
\end{equation}
Let $\vq\in X_0$ be a configuration. 
Let $\vY_j$ denote the $j$-th column of the  matrix $Y$. 
Then, by \eqref{eq:unknown}, 
$\vY_k$ belongs to the half-space of $X_0$ 
determined by the inequality $\sum_{i=1}^{n-1} x_i >0 \iff q_1-q_n>0$. 

As a consequence, for each $k=1,\ldots, n$ the vector
\[
\dfrac{\vY_j}
{Q_{1j}+Q_{jn}}
\]
belongs to $X_1$ (the sum of its components is $1$). 

\begin{theo}
Let $\psi\from \Delta^{n-2} \multimap X_1$ be the multi-valued map defined as follows:
for each $\vx \in \Delta^{n-2}$ the image $\psi(\vx) \subset X_1$ 
is the convex hull (which is a finite union of simplices) of the $n$ points
\[
\dfrac{\vY_j}
{Q_{1j}+Q_{jn}}
\]
for $j=1,\ldots, n$. 
Then
the map $\psi$ is continuous and 
there is a solution to the inverse central configuration problem 
for the configuration $\vx\in \Delta^{n-2}$ if and only if $\vx \in \psi(\vx)$. 
\end{theo}

\begin{remark}
\emph{A priori} the dimension of the simplices in $\psi(\vx)$ could be $<n-2$; 
as a corollary of Theorems (2.13) and (2.17) of \cite{ferrarioPfaffiansInverseProblem2020}, %
at least for $n$ small, 
the columns of $Q$ are always in general position (that is, the dimension
of the $n$ simplices of $\psi(\vx)$ is equal to $n-2$).
\end{remark}

\begin{example}[$n=3$]
For $n=3$, it is well-known that for $\alpha=1$ the inverse problem has solutions
for all configurations. It is in fact true for every $\alpha>0$, and 
  it can be seen for $\alpha=1$ in figure \ref{fig:2}.
\begin{figure}
\centering
\includegraphics[width=0.4\textwidth]{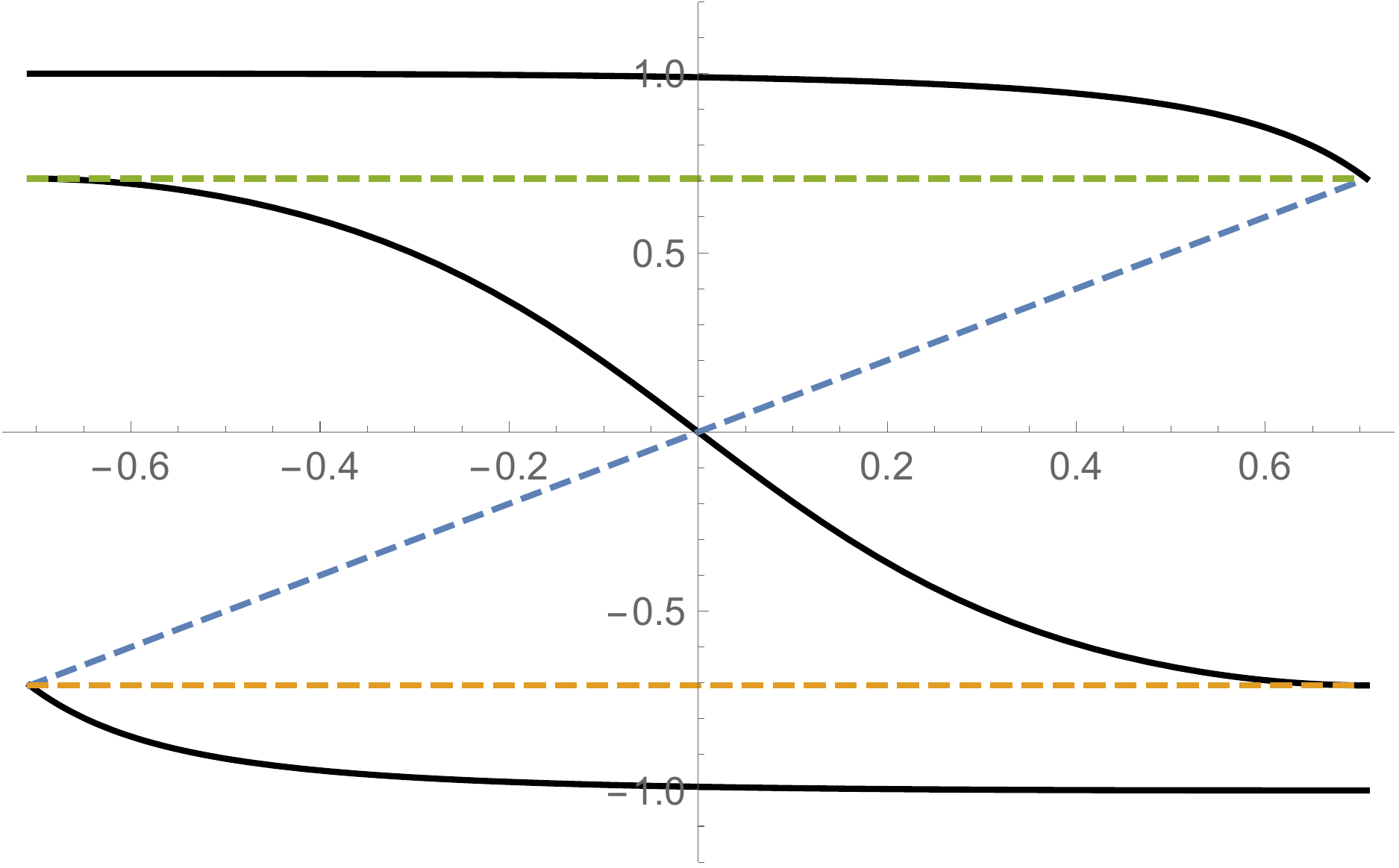}
\caption{The multi-valued map $\psi$ for $n=3$}
\label{fig:2}
\end{figure}
In order to visualize it, consider the homeomorphism $t \mapsto (t,\sqrt{1-t^2})$, which 
sends the interval 
$[-1,1]$ to the upper unit semicircle,
and the interval
$[-\dfrac{\sqrt{2}}{2}, \dfrac{\sqrt{2}}{2}]$ to the arc of the unit circle
with endpoints in $(\pm \dfrac{\sqrt{2}}{2}, \dfrac{\sqrt{2}}{2})$. Rotate in $\pi/4$ clockwise
and project it radially onto $\Delta^{1}$ (this last step is not strictly necessary): then  
consider the corresponding three columns $\vY_1$, $\vY_2$ and $\vY_3$ of the matrix $Y$, 
and their images in $X_1$, and the projections on the unit circle (actually, the
semicircle with $x_1+x_2>0$). Then rotating counter-clockwise by $\pi/4$ and 
taking the inverse of the homeomorphism $t \mapsto (t,\sqrt{1-t^2})$ gives 
the three maps 
\[
f_1, f_2, f_3 \from [-\dfrac{\sqrt{2}}{2}, \dfrac{\sqrt{2}}{2}]
\to [-1,1]
\]
represented in figure \ref{fig:2}, with $f_1>f_2 > f_3$. 
Now, it is clear that $\psi(\vx) = \CH[f_1(\vx),f_2(\vx)]$ 
(where $\CH$ means the convex hull of the following list of points), 
and if $t<0$ (and hence $x_1<x_2$) $\vx \in \CH[f_2(\vx), f_3(\vx)]$, 
while if $t>0$ (and hence $x_1>x_2$) $\vx \in \CH[f_1(\vx),f_2(\vx)]$. 
\end{example}

\begin{example}[$n=4$]
In this case the matrix \eqref{eq:matrixY} is
\[
\begin{bmatrix}
Q_{11} - Q_{21} & Q_{12} - Q_{22} & Q_{13}-Q_{23}  & Q_{14}-Q_{24} \\
Q_{21} - Q_{31} & Q_{22} - Q_{32} & Q_{23} - Q_{33} & Q_{24}-Q_{34} \\
Q_{31} - Q_{41} & Q_{32} - Q_{42} & Q_{33}-Q_{43} & Q_{34}-Q_{44} \\
\end{bmatrix}
\]
The $2$-simplex $\Delta^{n-2}$ and the plane $X_1$ can be projected on the unit
sphere, as in figure 
\ref{fig:3}. 
\begin{figure}
\centering
\includegraphics[width=0.4\textwidth]{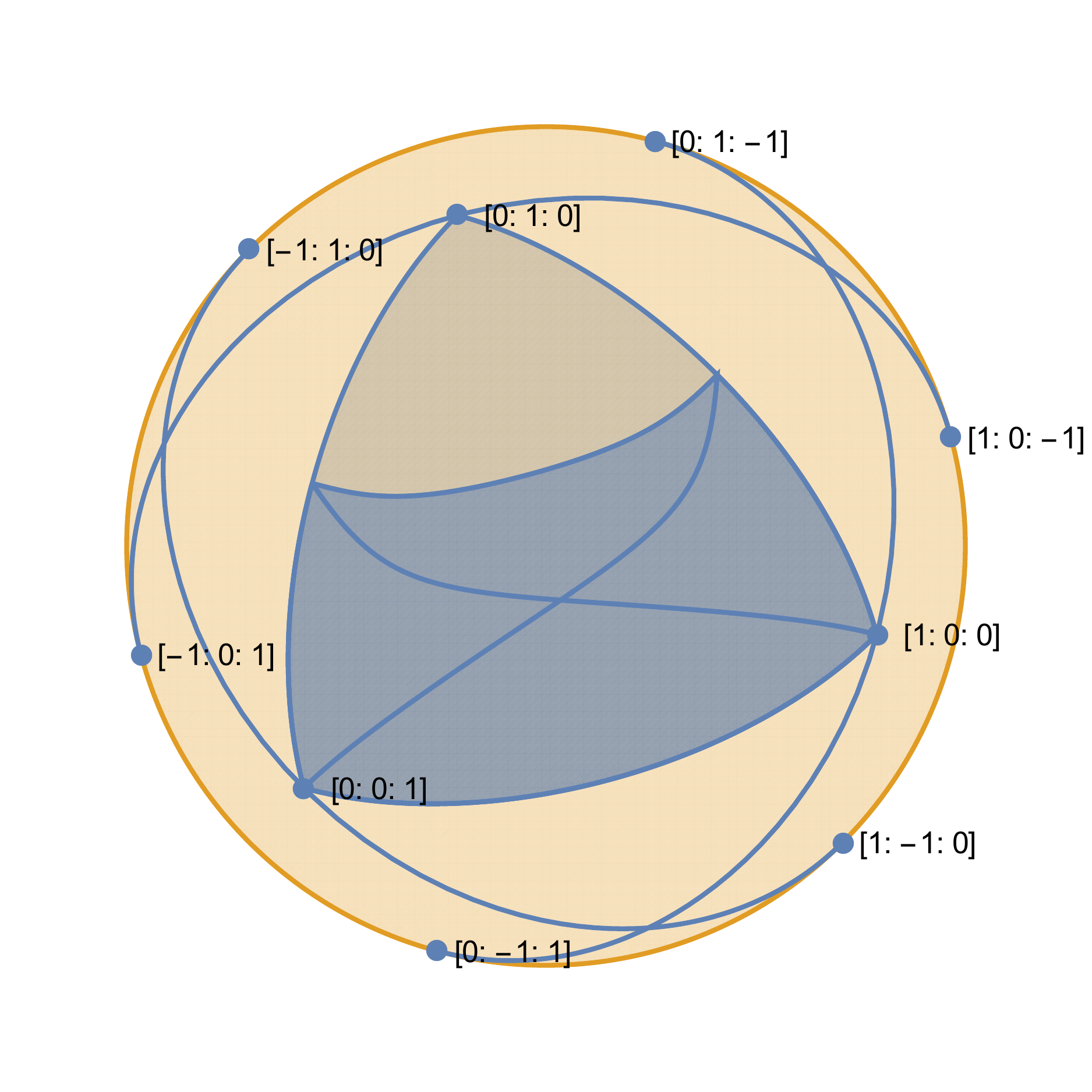}
\caption{The simplex $\Delta^2$ in $X_1$, for $n=4$ }
\label{fig:3}
\end{figure}
The simplex is represented together with the lines of equation $x_i+x_j=0$ in $X_1$. 
Now it is not possible to visualize the $4$ maps $f_1$, $f_2$, $f_3$ and $f_4$ from $\Delta^{2}$ 
to $X_1$,
defined as in the previous example. But after renormalization 
it is possible to consider the restrictions of $\psi$ to the three faces $x_i=0$ of $\Delta^2$, for $i=1,2,3$,
as follows.
The matrix as $x_1\to 0$ will be (up to norms of the columns):
\[
\begin{bmatrix}
1  & 1  & 0   & 0  \\
-1 & 0 & Q_{23}  & Q_{24}-Q_{34} \\
 0  & 0 & Q_{34} & Q_{34} \\
\end{bmatrix}
\]
since $Q_{1j} = Q_{2j}$ for each $j=3,4$. 
Recall that the matrix for the three particles $q_2,q_3,q_4$ is 
\[
\begin{bmatrix}
 Q_{22} - Q_{32} & Q_{23} - Q_{33} & Q_{24}-Q_{34} \\
 Q_{32} - Q_{42} & Q_{33}-Q_{43} & Q_{34}-Q_{44} \\
\end{bmatrix}
=
\begin{bmatrix}
 Q_{23} & Q_{23}  & Q_{24}-Q_{34} \\
Q_{24} -  Q_{23}  & Q_{34} & Q_{34} \\
\end{bmatrix}
\]
hence the limiting $Y$ as $x_1\to 0$ has a submatrix which corresponds
to the columns $\vY_2, \vY_3$ of the problem with $n=3$. 
The same computation can be performed for the other columns: 
it happens that the restrictions of $\psi$ to the faces of the simplex 
are nothing but the maps defined with $f_1$, $f_2$, and $f_3$ for suitable
choice of masses. For further details on such self-map,
see 
\cite{ferrarioPfaffiansInverseProblem2020}.
\end{example}

\subsection*{Acknowledgements}

I would like to thank the organizers and the participants 
of the \emph{Nielsen Theory and Related Topics}
conference 
held in 
Kortrijk 
June 3--7, 2019. In particular, I would like to thank 
Jan Andres, Peter Wong and 
Daciberg Gonçalves for 
very interesting discussions on this topic, among others. 

\noindent 
Address: 
DL Ferrario \\
Department of Mathematics and Applications\\
University of Milano--Bicocca\\
Via R. Cozzi 55\\
20125 Milano  -- Italy \\
email: \texttt{davide.ferrario@unimib.it}

\end{document}